\documentclass[12pt]{amsart}
\usepackage[margin=1in]{geometry}
\usepackage{amssymb,amsmath,amsthm,mathtools,amsfonts,times,hyperref,mathrsfs,multirow}
\usepackage{pgfplots}
\pgfplotsset{compat=1.15}
\usepackage{mathrsfs,mathdots}
\usetikzlibrary{arrows}
\usepackage{enumerate}
\usepackage{graphicx}
\usepackage{array}
\usepackage{ytableau}
\usepackage{pstricks,pst-node,graphicx}
\usepackage{tikz,varwidth}
\usepackage{setspace}
\usepackage{float}
\usetikzlibrary{calc}
\usepackage[super]{nth}
\usepackage{tikz}
\usepackage{extarrows}
\usepackage{pgfplots}
\pgfplotsset{compat=1.15}
\usetikzlibrary{arrows}

\newtheorem{theorem}{Theorem}[section]

\newtheorem{cor}[theorem]{Corollary}
\newtheorem{conj}[theorem]{Conjecture}

\theoremstyle{definition}

\theoremstyle{remark}
\newtheorem{remark}{Remark}
\numberwithin{equation}{section}

\allowdisplaybreaks

\newcommand\numberthis{\addtocounter{equation}{1}\tag{\theequation}}

\newcommand\restr[2]{{
  \left.\kern-\nulldelimiterspace 
  #1 
  \littletaller 
  \right|_{#2} 
  }}
\newcommand{\littletaller}{\mathchoice{\vphantom{\big|}}{}{}{}}  

\begin{document}
\title[]{Further Applications of Cubic $q$-Binomial Transformations}
\author{Alexander Berkovich}
\address{Department of Mathematics, University of Florida, Gainesville,
FL 32611, USA}
\email{alexb@ufl.edu}
\author{Aritram Dhar}
\address{Department of Mathematics, University of Florida, Gainesville,
FL 32611, USA}
\email{aritramdhar@ufl.edu}

\date{\today}



\begin{abstract}
Consider
\begin{align*}
G(N,M;\alpha,\beta,K,q) = \sum\limits_{j\in\mathbb{Z}}(-1)^jq^{\frac{1}{2}Kj((\alpha+\beta)j+\alpha-\beta)}\left[\begin{matrix}M+N\\N-Kj\end{matrix}\right]_{q}.    
\end{align*}
In this paper, we prove the non-negativity of coefficients of some cases of $G(N,M;\alpha,\beta,K,q)$. For instance, for non-negative integers $n$ and $t$, we prove that\\
\begin{align*}
G\left(n,n;\frac{4}{3}+\frac{3(3^t-1)}{2},\frac{5}{3}+\frac{3(3^t-1)}{2},3^{t+1},q\right)
\end{align*}
and
\begin{align*}
G\left(n-\frac{3^t-1}{2},n+\frac{3^t+1}{2};\frac{8}{3}+2(3^t-1),\frac{4}{3}-(3^t-1),3^{t+1},q\right)\\
\end{align*}
are polynomials in $q$ with non-negative coefficients. Using cubic positivity preserving transformations of Berkovich and Warnaar and some known formulae arising from Rogers-Szeg\"{o} polynomials, we establish new identities such as\\
\begin{align*}
\sum\limits_{0\le 3j\le n}\dfrac{(q^3;q^3)_{n-j-1}(1-q^{2n})q^{3j^2}}{(q;q)_{n-3j}(q^6;q^6)_{j}} = \sum\limits_{j=-\infty}^{\infty}(-1)^jq^{6j^2}{2n\brack n-3j}_q.
\end{align*}
\end{abstract}

\maketitle

\section{Introduction}\label{s1}
Let $L,m,n$ be non-negative integers. Define the conventional $q$-Pochhammer symbol as\\
\begin{align*}
(a)_L = (a;q)_L &:= \prod_{k=0}^{L-1}(1-aq^k),\\
(a)_{\infty} = (a;q)_{\infty} &:= \lim_{L\rightarrow \infty}(a)_L\,\,\text{where}\,\,\lvert q\rvert<1.\\
\end{align*}
Next, we define the $q$-binomial coefficient as\\
\begin{align*}
\left[\begin{matrix}m+n\\n\end{matrix}\right]_q := \Bigg\{\begin{array}{lr}
\dfrac{(q)_{m+n}}{(q)_m(q)_n}\quad\text{for } m,n\ge 0,\\
0\qquad\qquad\quad\text{otherwise}.\end{array}\\
\end{align*}
It is well-known that $\left[\begin{matrix}m+n\\n\end{matrix}\right]_q$ is the generating function for partitions into at most $n$ parts each of size at most $m$ (see \cite{A98}).\\\par Throughout the remainder of the paper, $P(q)\ge 0$ means that a power series in $q$, $P(q)$, has non-negative coefficients.\\\par For non-negative integers $N, M$, positive integers $i, K$ such that $i < K$, and $\alpha, \beta\ge 0$, define\\
\begin{align*}
D_{K,i}&(N,M;\alpha,\beta;q) = D_{K,i}(N,M;\alpha,\beta)\\
&:= \sum\limits_{j\in\mathbb{Z}}\Bigg\{q^{j((\alpha+\beta)Kj+K\beta-(\alpha+\beta)i)}\left[\begin{matrix}M+N\\M-Kj\end{matrix}\right]_{q}-q^{((\alpha+\beta)j+\beta)(Kj+i)}\left[\begin{matrix}M+N\\M-Kj-i\end{matrix}\right]_{q}\Bigg\}.\numberthis\label{eq11}\\
\end{align*}
Andrews, Baxter, Bressoud, Burge, Forrester, Viennot \cite{ABBBFV87} showed that $D_{K,i}(N,M;\alpha,\beta)$ is the generating function for a certain class of restricted partitions when $\alpha,\beta\in\mathbb{N}\cup\{0\}$, $1\le \alpha+\beta\le K-1$, and $\beta-i\le N-M\le K-\alpha-i$. Thus,\\
\begin{align*}
D_{K,i}(N,M;\alpha,\beta)\ge 0.\numberthis\label{eq12}\\
\end{align*}
\par Bressoud \cite{B96} then considered the following polynomials\\
\begin{align*}
G(N,M;\alpha,\beta,K,q) = G(N,M;\alpha,\beta,K) &:= D_{2K,K}(N,M;\alpha,\beta)\\
&= \sum\limits_{j\in\mathbb{Z}}(-1)^jq^{\frac{1}{2}Kj((\alpha+\beta)j+\alpha-\beta)}\left[\begin{matrix}M+N\\N-Kj\end{matrix}\right]_{q}.\numberthis\label{eq13}\\
\end{align*}
and made the following conjecture \cite[Conjecture $6$]{B96}.\\
\begin{conj}\label{conj11}
Let $K$ be a positive integer and $N,M,\alpha K,\beta K$ be non-negative integers such that $1\le \alpha+\beta\le 2K-1$ (strict inequalities when $K = 2$) and $\beta-K\le N-M\le K-\alpha$. Then, $G(N,M;\alpha,\beta,K)$ is a polynomial in $q$ with non-negative coefficients.\\
\end{conj}
Many cases of Conjecture \ref{conj11} were proven in the literature \cite{B20,B22,BD25,BW05,B81,IKS99,W22,W01,W03}.\\\par Note that one of the mod $3$ conjectures due to Borwein \cite{A95} can be stated as\\
\begin{equation*}
A_n(q) := G(n,n;4/3,5/3,3)\ge 0,\numberthis\label{eq14}\\
\end{equation*}
\begin{equation*}
B_n(q) := G(n+1,n-1;2/3,7/3,3)\ge 0,\numberthis\label{eq15}\\
\end{equation*}
and
\begin{equation*}
C_n(q) := G(n+1,n-1;1/3,8/3,3)\ge 0.\numberthis\label{eq16}
\end{equation*}\\
All the three inequalities above were proven by Wang \cite{W22} and Wang and Krattenthaler \cite{KW22}.\\\par In 2020, Berkovich \cite{B20} showed that\\
\begin{equation*}
G(n,n+1;8/3,4/3,3) = \sum\limits_{k=0}^{n}q^{\frac{k(k+1)}{2}}\left[\begin{matrix}n\\k\end{matrix}\right]_{q}(-q)_k\ge 0,\numberthis\label{eq17}\\    
\end{equation*}
and
\begin{equation*}
G(n,n+1;4/3,2/3,3) = \sum\limits_{k=0}^{n}q^{(n+1)k}\left[\begin{matrix}n\\k\end{matrix}\right]_{q}(-q)_{n-k}\ge 0.\numberthis\label{eq18}
\end{equation*}\\
Note that \eqref{eq18} follows from \eqref{eq17} using the transformation $q\rightarrow q^{-1}$.\\\quad\\
Recently, Berkovich and Dhar \cite{BD25} gave the following generalized conjecture regarding non-negativity of $D_{K,i}(N,M;\alpha,\beta)$.\\
\begin{conj}\label{conj12}
Let $K,i$ be positive integers such that $0 < i < K$ and $N,M,\alpha K,\beta K,\alpha i,\beta i$ be non-negative integers such that $1\le \alpha+\beta\le K-1$ (strict inequalities when $K = 4$ and $i=2$) and $\beta-i\le N-M\le K-\alpha-i$. Then, $D_{K,i}(N,M;\alpha,\beta)$ is a polynomial in $q$ with non-negative coefficients.\\
\end{conj}
It is easy to see that Conjecture \ref{conj11} is the special case $(i,K)\mapsto (K,2K)$ of Conjecture \ref{conj12}.\\\par Berkovich and Dhar proved some special cases of Conjecture \ref{conj12} in \cite{BD25} using certain \textit{positivity-preserving transformations} for $q$-binomial coefficients due to Berkovich and Warnaar \cite{BW05}. In particular, we will focus our attention on the following two \textit{cubic positivity-preserving transformations} from \cite{BW05}.\\
\begin{theorem}(\cite[Lemma $2.6$]{BW05} $L$, $j$, $r$ even case)\label{thm13}
For integers $L$ and $j$, we have\\
\begin{equation}\label{eq19}
\sum\limits_{r=0}^{\left\lfloor\frac{L}{3}\right\rfloor}T_{L,r}(q)\left[\begin{matrix}2r \\ r-j\end{matrix}\right]_{q^3} = q^{3j^2}\left[\begin{matrix}2L \\ L-3j\end{matrix}\right]_{q},
\end{equation}\\
where\\
\begin{equation}\label{eq110}
T_{L,r}(q) = \dfrac{q^{3r^2}(q^3;q^3)_{L-r-1}(1-q^{2L})}{(q^3;q^3)_{2r}(q;q)_{L-3r}}.
\end{equation}\\
When $L=r=0$, $T_{0,0}(q) := 1$.\\
\end{theorem}
Berkovich and Warnaar \cite{BW05} showed that\\
\begin{equation}\label{eq111}
f_{L,r}(q) = \dfrac{(q^3;q^3)_{\frac{1}{2}(L-r-2)}(1-q^L)}{(q^3;q^3)_r(q;q)_{\frac{1}{2}(L-3r)}}
\end{equation}\\ 
is a  polynomial with non-negative coefficients for $0\le 3r\le L$ and $r\equiv L\pmod 2$. It is then evident from \eqref{eq19} and \eqref{eq110} that\\ $$T_{L,r}(q) = q^{3r^2}f_{2L,2r}(q)\\$$ has non-negative coefficients.\\\par
It is then easy to verify that for any identity of the form\\
\begin{equation}\label{eq112}
F_T(L,q) = \sum\limits_{j\in\mathbb{Z}}\alpha(j,q)\left[\begin{matrix}2L \\ L-j\end{matrix}\right]_{q^3},    
\end{equation}\\
using transformation \eqref{eq19}, the following identity holds\\
\begin{equation}\label{eq113}
\sum\limits_{r\ge 0}T_{L,r}(q)F_T(r,q) = \sum\limits_{j\in\mathbb{Z}}\alpha(j,q)\sum\limits_{r\ge 0}T_{L,r}(q)\left[\begin{matrix}2r \\ r-j\end{matrix}\right]_{q^3} = \sum\limits_{j\in\mathbb{Z}}\alpha(j,q)q^{3j^2}\left[\begin{matrix}2L \\ L-3j\end{matrix}\right]_{q}.
\end{equation}\\
Hence, if $F_T(L,q)\ge 0$, then\\
\begin{equation}\label{eq114}
\sum\limits_{j\in\mathbb{Z}}\alpha(j,q)q^{3j^2}\left[\begin{matrix}2L \\ L-3j\end{matrix}\right]_{q}\ge 0.    
\end{equation}\\
So, we say that transformation \eqref{eq19} is positivity-preserving.\\
\begin{theorem}(\cite[Lemma $2.6$]{BW05} $L$, $j$, $r$ odd case)\label{thm14}
For integers $L$ and $j$, we have\\
\begin{equation}\label{eq115}
\sum\limits_{r=0}^{\left\lfloor\frac{L}{3}\right\rfloor}\Tilde{T}_{L,r}(q)\left[\begin{matrix}2r+1 \\ r-j\end{matrix}\right]_{q^3} = q^{3j^2+3j}\left[\begin{matrix}2L+1 \\ L-3j-1\end{matrix}\right]_{q},
\end{equation}\\
where\\
\begin{equation}\label{eq116}
\Tilde{T}_{L,r}(q) = \dfrac{q^{3r^2+3r}(q^3;q^3)_{L-r-1}(1-q^{2L+1})}{(q^3;q^3)_{2r+1}(q;q)_{L-3r-1}}.
\end{equation}\\
When $L=r=0$, $\Tilde{T}_{0,0}(q) := 0$.\\
\end{theorem}
It is then evident from \eqref{eq111} and \eqref{eq115} that\\ $$\Tilde{T}_{L,r}(q) = q^{3r^2+3r}f_{2L+1,2r+1}(q)\\$$ has non-negative coefficients.\\\par
It is then easy to verify that for any identity of the form\\
\begin{equation}\label{eq117}
F_{\Tilde{T}}(L,q) = \sum\limits_{j\in\mathbb{Z}}\alpha(j,q)\left[\begin{matrix}2L+1 \\ L-j\end{matrix}\right]_{q^3},    
\end{equation}\\
using transformation \eqref{eq114}, the following identity holds\\
\begin{equation}\label{eq118}
\sum\limits_{r\ge 0}\Tilde{T}_{L,r}(q)F_{\Tilde{T}}(r,q) = \sum\limits_{j\in\mathbb{Z}}\alpha(j,q)\sum\limits_{r\ge 0}\Tilde{T}_{L,r}(q)\left[\begin{matrix}2r+1 \\ r-j\end{matrix}\right]_{q^3} = \sum\limits_{j\in\mathbb{Z}}\alpha(j,q)q^{3j^2+3j}\left[\begin{matrix}2L+1 \\ L-3j-1\end{matrix}\right]_{q}.    
\end{equation}\\
Hence, if $F_{\Tilde{T}}(L,q)\ge 0$, then\\
\begin{equation}\label{eq119}
\sum\limits_{j\in\mathbb{Z}}\alpha(j,q)q^{3j^2+3j}\left[\begin{matrix}2L+1 \\ L-3j-1\end{matrix}\right]_{q}\ge 0. 
\end{equation}\\
Again, we say that transformation \eqref{eq114} is positivity-preserving.\\\par In an attempt to prove Borwein's mod $3$ conjecture, Andrews \cite[Theorem $4.1$]{A95} gave the following identities.\\
\begin{theorem}\label{thm15}
For $n > 0$, we have\\
\begin{align*}
A_n(q) &= \sum\limits_{0\le 3j\le n}\dfrac{(q^3;q^3)_{n-j-1}(1-q^{2n})(q;q)_{3j}q^{3j^2}}{(q;q)_{n-3j}(q^3;q^3)_{2j}(q^3;q^3)_j},\numberthis\label{eq120}\\
B_n(q) &= \sum\limits_{0\le 3j\le n-1}\dfrac{(q^3;q^3)_{n-j-1}(1-q^{3j+2}-q^{n+3j+2}+q^{n+1})(q;q)_{3j}q^{3j^2+3j}}{(q;q)_{n-3j-1}(q^3;q^3)_{2j+1}(q^3;q^3)_j},\numberthis\label{eq121}\\
C_n(q) &= \sum\limits_{0\le 3j\le n-1}\dfrac{(q^3;q^3)_{n-j-1}(1-q^{3j+1}-q^{n+3j+2}+q^n)(q;q)_{3j}q^{3j^2+3j}}{(q;q)_{n-3j-1}(q^3;q^3)_{2j+1}(q^3;q^3)_j}.\numberthis\label{eq122}\\
\end{align*}
where $A_n(q)$, $B_n(q)$, and $C_n(q)$ are defined in \eqref{eq14}, \eqref{eq15}, and \eqref{eq16} respectively.\\
\end{theorem}
From \eqref{eq120}, \eqref{eq121}, and \eqref{eq122} above, it is not clear that $A_n(q)$, $B_n(q)$, and $C_n(q)$ are non-negative.\\\par Now, we state new identities which are similar to the identities in Theorem \ref{thm15}.\\
\begin{theorem}\label{thm16}
For $n > 0$ and $a\in\{0,1\}$, we have\\
\begin{align*}
\sum\limits_{0\le 3j\le n-a}\dfrac{(-1)^j(q^3;q^3)_{n-j-1}(1-q^{2n+a})q^{3j^2}}{(q;q)_{n-3j-a}(q^6;q^6)_{j}} = \sum\limits_{j=-\infty}^{\infty}(-1)^jq^{3j^2}{2n+a\brack n-3j-a}_q.\numberthis\label{eq123}\\
\end{align*}
\end{theorem}
\begin{theorem}\label{thm17}
For $n > 0$, we have\\
\begin{align*}
\sum\limits_{0\le 3j\le n}\dfrac{(-1)^j(q^3;q^3)_{n-j-1}(1-q^{2n})q^{3j^2-3j}}{(q;q)_{n-3j}(q^6;q^6)_{j}} = \sum\limits_{j=-\infty}^{\infty}(-1)^jq^{3j^2+3j}{2n\brack n-3j}_q.\numberthis\label{eq124}\\
\end{align*}
\end{theorem}
\begin{theorem}\label{thm18}
For $n > 0$, we have\\
\begin{align*}
\sum\limits_{0\le 3j\le n}\dfrac{(q^3;q^3)_{n-j-1}(1-q^{2n})q^{3j^2}}{(q;q)_{n-3j}(q^6;q^6)_{j}} = \sum\limits_{j=-\infty}^{\infty}(-1)^jq^{6j^2}{2n\brack n-3j}_q.\numberthis\label{eq125}\\
\end{align*}
\end{theorem}
\begin{theorem}\label{thm19}
For $n > 0$ and $a\in\{0,1\}$, we have\\
\begin{align*}
\sum\limits_{0\le 3j\le n-a}\dfrac{(q^3;q^3)_{n-j-1}(1-q^{2n+a})q^{3j^2+3j}}{(q;q)_{n-3j-a}(q^6;q^6)_{j}} = \sum\limits_{j=-\infty}^{\infty}(-1)^jq^{6j^2+3j}{2n+a\brack n-3j-a}_q.\numberthis\label{eq126}\\
\end{align*}
\end{theorem}
\begin{remark}\label{rmk1}
It is to be noted here that the right-hand sides of \eqref{eq123}-\eqref{eq126} are non-negative. These follow from \eqref{eq12}. However, the left-hand sides of \eqref{eq113}-\eqref{eq126} are not obvious to be non-negative. The new sums on the left-hand sides are different compared to Equations $(3.4)$, $(3.5)$, and $(3.6)$ in Andrews' paper \cite{A95} as they contain the factor $1/(q^6;q^6)_j$ in the summand (among other differences) and the main difference between the formal bilateral sums on the right-hand sides is the quadratic power of $q$ appearing in them. The bilateral sums in Theorems \ref{thm16}-\ref{thm19} contain $q^{3j^2}$ or $q^{6j^2}$ as a factor, while those in Andrews' paper \cite{A95} contain $q^{\frac{9}{2}j^2}$. (Interestingly, $q^{\frac{9}{2}j^2}$ is exactly half way between $3j^2$ and $6j^2$.)\\    
\end{remark}
We now state an identity which is different from those in Theorems \ref{thm16}-\ref{thm19} in the sense the right-hand side does not have a $(-1)^j$ factor in the summand and hence, is manifestly positive whereas the right-hand sides of Theorems \ref{thm16}-\ref{thm19} are $D_{K,i}(N,M;\alpha,\beta)$ functions where $D_{K,i}(N,M;\alpha,\beta)$ are defined in \eqref{eq11}. It is as follows. 
\begin{theorem}\label{thm110}
For $n > 0$ and $a\in\{0,1\}$, we have\\
\begin{align*}
\sum\limits_{0\le 3j\le n-a}\dfrac{(q^6;q^6)_{n-j-1}(-q^3;q^3)_{2j+a}(1-q^{4n+2a})q^{6j^2+(6a-3)j-3a}}{(q^2;q^2)_{n-3j-a}(q^6;q^6)_{2j+a}} = \sum\limits_{j=-\infty}^{\infty}q^{6j^2+(6a+3)j}{2n+a\brack n-3j-a}_{q^2}.\numberthis\label{eq127}\\
\end{align*}    
\end{theorem}
Also note that the $q$-binomial coefficient in \eqref{eq127} has base $q^2$ instead of base $q$ as in Theorems \ref{thm16}-\ref{thm19}.\\\par We now state two general inequalities.\\
\begin{theorem}\label{thm111}
For non-negative integers $n$, $t$, $x$, $y$ and any integer $a$,\\
\begin{align*}
G\left(n+3^ta,n-3^ta;\frac{x}{3}+\frac{(3^t-1)(3-2a)}{2},\frac{y}{3}+\frac{(3^t-1)(3+2a)}{2},3^{t+1}\right)\ge 0\numberthis\label{eq128}     
\end{align*}
if $G(n+a,n-a;x/3,y/3,3)\ge 0$.\\
\end{theorem}
\begin{theorem}\label{thm112}
For non-negative integers $n$, $t$, $x$, $y$ and any integer $a$,\\
\begin{align*}
G\left(n-3^ta-\frac{3^t-1}{2},n+3^ta+\frac{3^t+1}{2};\frac{x}{3}+(3^t-1)(a+2),\frac{y}{3}+(3^t-1)(a-1),3^{t+1}\right)\ge 0\numberthis\label{eq129}    
\end{align*}
if $G(n-a,n+a+1;x/3,y/3,3)\ge 0$.\\
\end{theorem}
We conclude this section with the following important corollaries.\\
\begin{cor}\label{cor113}
For non-negative integers $n$ and $t$, we have\\
\begin{align*}
G\left(n,n;\frac{4}{3}+\frac{3(3^t-1)}{2},\frac{5}{3}+\frac{3(3^t-1)}{2},3^{t+1}\right)\ge 0.\numberthis\label{eq130}\\    
\end{align*}
\begin{align*}
G\left(n+3^t,n-3^t;\frac{2}{3}+\frac{3^t-1}{2},\frac{7}{3}+\frac{5(3^t-1)}{2},3^{t+1}\right)\ge 0.\numberthis\label{eq131}\\ 
\end{align*}
\begin{align*}
G\left(n+3^t,n-3^t;\frac{1}{3}+\frac{3^t-1}{2},\frac{8}{3}+\frac{5(3^t-1)}{2},3^{t+1}\right)\ge 0.\numberthis\label{eq132}\\ 
\end{align*}
\end{cor}
\begin{cor}\label{cor114}
For non-negative integers $n$ and $t$, we have\\
\begin{align*}
G\left(n-\frac{3^t-1}{2},n+\frac{3^t+1}{2};\frac{8}{3}+2(3^t-1),\frac{4}{3}-(3^t-1),3^{t+1}\right)\ge 0.\numberthis\label{eq133}\\    
\end{align*}
\begin{align*}
G\left(n-\frac{3^t-1}{2},n+\frac{3^t+1}{2};\frac{4}{3}+2(3^t-1),\frac{2}{3}-(3^t-1),3^{t+1}\right)\ge 0.\numberthis\label{eq134}\\    
\end{align*}
\end{cor}

\section{Proofs}\label{s2}
In this section, we provide proofs of our main results stated in \S\ref{s1}.\\
\subsection{Proofs of Theorems \ref{thm16}-\ref{thm110}}
We start by defining the \textit{Rogers-Szeg\"{o} polynomials}. For any non-negative integer $n$, the \textit{Rogers-Szeg\"{o} polynomials} are defined as \cite[Ch. $3$, Examples $3$-$9$]{A98}\\
\begin{align*}
H_n(t;q) = H_n(t) := \sum\limits_{j=0}^{n}t^j\left[\begin{matrix}n\\j\end{matrix}\right]_{q}.\numberthis\label{eq21}\\    
\end{align*}
Then the following special cases are well-known \cite{A98,BW05}.\\
\begin{align*}
H_{2n}(-1) = (q;q^2)_n,\numberthis\label{eq22}    
\end{align*}
and
\begin{align*}
H_{n}(-q) = (q;q^2)_{\lfloor (n+1)/2\rfloor}.\numberthis\label{eq23}\\    
\end{align*}
Another well-known special case of the \textit{Rogers-Szeg\"{o} polynomials} is the evaluation \cite[p. $49$, Chapter $3$, Example $5$]{A98}\\
\begin{align*}
H_n(q^{\frac{1}{2}}) = (-q^{\frac{1}{2}};q^{\frac{1}{2}})_n.\numberthis\label{eq24}\\
\end{align*}
It is easy to show that \eqref{eq22} can be re-written as\\
\begin{align*}
\sum\limits_{j=-n}^{n}(-1)^j\left[\begin{matrix}2n\\n-j\end{matrix}\right]_{q} = (-1)^n(q;q^2)_n.\numberthis\label{eq25}\\    
\end{align*}
Now, substituting $q\mapsto q^3$ in \eqref{eq25} and applying \eqref{eq113}, we get \eqref{eq123} with $a=0$.\\\par Replacing $n\mapsto 2n+1$ in \eqref{eq23}, we can re-write \eqref{eq23} as\\
\begin{align*}
\sum\limits_{j=-n-1}^{n}(-1)^jq^j\left[\begin{matrix}2n+1\\n-j\end{matrix}\right]_{q} = (-1)^{n+1}q^{-n-1}(q;q^2)_{n+1}.\numberthis\label{eq26}\\    
\end{align*}
Then, substituting $q\mapsto q^3$ in \eqref{eq26} and applying \eqref{eq118}, we get \eqref{eq123} with $a=1$ which completes the proof of Theorem \ref{thm16}.\\\par Similarly, replacing $n\mapsto 2n$ in \eqref{eq23}, we can re-write \eqref{eq23} as\\
\begin{align*}
\sum\limits_{j=-n}^{n}(-1)^jq^j\left[\begin{matrix}2n\\n-j\end{matrix}\right]_{q} = (-1)^nq^{-n}(q;q^2)_n.\numberthis\label{eq27}\\    
\end{align*}
Then, substituting $q\mapsto q^3$ in \eqref{eq27} and applying \eqref{eq113}, we get \eqref{eq124} which proves Theorem \ref{thm17}.\\\par Now, replacing $q\mapsto q^{-1}$ in \eqref{eq25}, we get\\
\begin{align*}
\sum\limits_{j=-n}^{n}(-1)^jq^{j^2}\left[\begin{matrix}2n\\n-j\end{matrix}\right]_{q} = (q;q^2)_n.\numberthis\label{eq28}\\    
\end{align*}
\eqref{eq28} was also obtained by Andrews in \cite[eq. $(2.2)$]{A91}. Now, substituting $q\mapsto q^3$ in \eqref{eq28} and applying \eqref{eq113}, we get \eqref{eq125} which proves Theorem \ref{thm18}.\\\par Replacing $q\mapsto q^{-1}$ in \eqref{eq26}, we get\\
\begin{align*}
\sum\limits_{j=-n-1}^{n}(-1)^jq^{j^2}\left[\begin{matrix}2n+1\\n-j\end{matrix}\right]_{q} = (q;q^2)_{n+1}.\numberthis\label{eq29}\\    
\end{align*}
Now, substituting $q\mapsto q^3$ in \eqref{eq29} and applying \eqref{eq118}, we get \eqref{eq126} with $a=1$.\\\par Similarly, replacing $q\mapsto q^{-1}$ in \eqref{eq27}, we get\\
\begin{align*}
\sum\limits_{j=-n}^{n}(-1)^jq^{j^2+j}\left[\begin{matrix}2n\\n-j\end{matrix}\right]_{q} = q^n(q;q^2)_n.\numberthis\label{eq210}\\    
\end{align*}
Now, substituting $q\mapsto q^3$ in \eqref{eq210} and applying \eqref{eq113}, we get \eqref{eq126} with $a=0$ which proves Theorem \ref{thm19}.\\\par Replacing $n\mapsto 2n$, it is easy to show that \eqref{eq24} can be re-written as\\
\begin{align*}
\sum\limits_{j=-n}^{n}q^{\frac{j}{2}}\left[\begin{matrix}2n\\n-j\end{matrix}\right]_{q} = q^{-\frac{n}{2}}(-q^{\frac{1}{2}};q^{\frac{1}{2}})_{2n}.\numberthis\label{eq211}\\    
\end{align*}
Now, substituting $q\mapsto q^3$ in \eqref{eq211} and applying \eqref{eq113}, we get \eqref{eq127} with $a=0$ after replacing $q$ by $q^2$.\\\par Similarly, replacing $n\mapsto 2n+1$, it is easy to show that \eqref{eq24} can be re-written as\\
\begin{align*}
\sum\limits_{j=-n-1}^{n}q^{\frac{j}{2}}\left[\begin{matrix}2n+1\\n-j\end{matrix}\right]_{q} = q^{\frac{-n-1}{2}}(-q^{\frac{1}{2}};q^{\frac{1}{2}})_{2n+1}.\numberthis\label{eq212}\\    
\end{align*}
Now, substituting $q\mapsto q^3$ in \eqref{eq212} and applying \eqref{eq118}, we get \eqref{eq127} with $a=1$ after replacing $q$ by $q^2$ which proves Theorem \ref{thm110}.\\\qed\\

\subsection{Proofs of Theorems \ref{thm111} \& \ref{thm112}}
We begin by assuming that\\
\begin{align*}
G\left(n+a,n-a;\frac{x}{3},\frac{y}{3},3,q\right) = \sum\limits_{j=-\infty}^{\infty}(-1)^jq^{\frac{(x+y)j^2+(x-y)j}{2}}\left[\begin{matrix}2n\\n+a-3j\end{matrix}\right]_{q}\ge 0,\numberthis\label{eq213}\\    
\end{align*}
where the conditions for non-negativity in \eqref{eq213} follow from those in Conjecture \ref{conj11}. Making the substitution $q\mapsto q^3$ in \eqref{eq213} and applying \eqref{eq113}, we get\\
\begin{align*}
\sum\limits_{r=0}^{\left\lfloor\frac{n}{3}\right\rfloor}T_{n,r}(q)G\left(r+a,r-a;\frac{x}{3},\frac{y}{3},3,q^3\right) = q^{3a^2}G\left(n+3a,n-3a;\frac{x}{3}+3-2a,\frac{y}{3}+3+2a,9,q\right).\numberthis\label{eq214}\\
\end{align*}
Since $T_{n,r}(q)\ge 0$, we have\\
\begin{align*}
G\left(n+3a,n-3a;\frac{x}{3}+3-2a,\frac{y}{3}+3+2a,9,q\right)\ge 0.\numberthis\label{eq215}\\
\end{align*}
Now, iterating the same process $t$ $(\ge 0)$ times, we get \eqref{eq128}.\\\par Similarly, we assume that\\
\begin{align*}
G\left(n-a,n+a+1;\frac{x}{3},\frac{y}{3},3,q\right) = \sum\limits_{j=-\infty}^{\infty}(-1)^jq^{\frac{(x+y)j^2+(x-y)j}{2}}\left[\begin{matrix}2n+1\\n-a-3j\end{matrix}\right]_{q}\ge 0,\numberthis\label{eq216}\\    
\end{align*}
where the conditions for non-negativity in \eqref{eq216} follow from those in Conjecture \ref{conj11}. Making the substitution $q\mapsto q^3$ in \eqref{eq216} and applying \eqref{eq118}, we get\\
\begin{equation}\label{eq217}
\begin{multlined}
\sum\limits_{r=0}^{\left\lfloor\frac{n}{3}\right\rfloor}\Tilde{T}_{n,r}(q)G\left(r-a,r+a+1;\frac{x}{3},\frac{y}{3},3,q^3\right)\\ = q^{3a^2+3a}G\left(n-3a-1,n+3a+2;\frac{x}{3}+2(a+2),\frac{y}{3}+2(a-1),9,q\right).\\
\end{multlined}
\end{equation}
Since $\Tilde{T}_{n,r}(q)\ge 0$, we have\\
\begin{align*}
G\left(n-3a-1,n+3a+2;\frac{x}{3}+2(a+2),\frac{y}{3}+2(a-1),9,q\right)\ge 0.\numberthis\label{eq218}\\
\end{align*}
Now, iterating the same process $t$ $(\ge 0)$ times, we get \eqref{eq129}.\\\qed\\

\subsection{Proofs of Corollaries \ref{cor113} \& \ref{cor114}}
\eqref{eq130} follows from \eqref{eq14} and the substitution $(a,x,y) = (0,4,5)$ in \eqref{eq128}. \eqref{eq131} follows from \eqref{eq15} and the substitution $(a,x,y) = (1,2,7)$ in \eqref{eq128}. \eqref{eq132} follows from \eqref{eq16} and the substitution $(a,x,y) = (1,1,8)$ in \eqref{eq128}. This completes the proof of Corollary \ref{cor113}.\\\par Similarly, \eqref{eq133} follows from \eqref{eq17} and the substitution $(a,x,y) = (0,8,4)$ in \eqref{eq129}. \eqref{eq134} follows from \eqref{eq18} and the substitution $(a,x,y) = (0,4,2)$ in \eqref{eq129}. This completes the proof of Corollary \ref{cor114}.\\\qed

\section*{Acknowledgments}\label{s4}
We would like to thank George E. Andrews for his kind interest and for directing us to \cite[eq. $(2.2)$]{A91}. We would also like to thank the anonymous referee for his/her comments and suggestions.

\bibliographystyle{amsplain}

\begin{thebibliography}{10}


\bibitem{A91}
G.~E.~Andrews,
Partitions and the Gaussian sum,
The mathematical heritage of C. F. Gauss,
World Sci. Publ.,
River Edge, NJ,
1991,
pp.35--42.


\bibitem{A95}
G.~E.~Andrews,
On a Conjecture of Peter Borwein,
\textit{J. Symbolic Comput.}
20 (5-6)
(1995)
487--501.


\bibitem{A98}
G.~E.~Andrews,
The Theory of Partitions,
Cambridge University Press,
1998.


\bibitem{ABBBFV87}
G.~E.~Andrews, R.~J.~Baxter, D.~M.~Bressoud, W.~H.~Burge, P.~J.~Forrester, G.~Viennot, Partitions with prescribed hook differences,
\textit{Eur. J. Comb.}
8 (4)
(1987)
341--350.


\bibitem{B20}
A.~Berkovich,
Some new positive observations,
\textit{Discrete Math.}
343 (11)
(2020) 112040,
8 pp.


\bibitem{B22}
A.~Berkovich,
Bressoud's identities for even moduli. New companions and related positivity results,
\textit{Discrete Math.}
345 (12)
(2022) 113104,
13 pp.


\bibitem{BD25}
A.~Berkovich, A.~Dhar,
Extension of Bressoud's generalization of Borwein's Conjecture and some exact results, \textit{Ramanujan J.}
67 (16),
2025
(to appear in the special issue in honor of the 85th birthdays of George E. Andrews and Bruce C. Berndt).


\bibitem{BW05}
A.~Berkovich, S.~O.~Warnaar,
Positivity preserving transformations for $q$-binomial coefficients,
\textit{Trans. Am. Math. Soc.}
\textbf{357 (6)}
(2005)
2291--2351.


\bibitem{B81}
D.~M.~Bressoud,
Some identities for terminating $q$-series,
\textit{Math. Proc. Camb. Philos. Soc.}
\textbf{89}
(1981)
211--223.


\bibitem{B96}
D.~M.~Bressoud,
The Borwein conjecture and partitions with prescribed hook differences,
\textit{Electron. J. Combin.}
\textbf{3 (2)}
(1996)
\#4.


\bibitem{IKS99}
M.~E.~H.~Ismail, D.~Kim, D.~Stanton,
Lattice paths and positive trigonometric sums,
\textit{Constr. Approx.}
\textbf{15}
(1999)
69--81.


\bibitem{W22}
C.~Wang,
Analytic proof of the Borwein conjecture,
\textit{Adv. Math.}
\textbf{394}
(2022) 108028
54 pp.


\bibitem{KW22}
C.~Wang, C.~Krattenthaler,
An asymptotic approach to Borwein-type sign pattern theorems,
preprint,
\href{https://arxiv.org/abs/2201.12415}{arXiv:2201.12415}.


\bibitem{W01}
S.~O.~Warnaar,
The generalized Borwein conjecture. I. The Burge transform,
in: B. C. Berndt, K. Ono (Eds.),
\textit{q-Series with Applications to Combinatorics, Number Theory and Physics,
in: Contemp. Math},
vol. 291, AMS, Providence, RI,
2001,
pp. 243--267.


\bibitem{W03}
S.~O.~Warnaar,
The generalized Borwein conjecture. II. Refined $q$-trinomial coefficients,
\textit{Discrete Math.}
\textbf{272 (2-3)}
(2003)
215--258.


\end{thebibliography}


\end{document}